\documentclass[12pt]{amsart}  
\usepackage{hyperref}

\usepackage{amssymb,amsmath,amsthm,amscd}
\usepackage[all]{xy}

\numberwithin{equation}{section}

\newtheoremstyle{fancy1}{10pt}{10pt}{\itshape}{12pt}{\textsc\bgroup}{.\egroup}{8pt}{ }
\newtheoremstyle{fancy2}{10pt}{10pt}{}{12pt}{\itshape}{.}{8pt}{ }

\theoremstyle{fancy1}

\newtheorem{lem}[equation]{Lemma}
\newtheorem{prop}[equation]{Proposition}

\newtheorem{main}{Theorem}
\newtheorem*{main*}{Theorem}

\newtheorem*{cor*}{Corollary}

\newtheorem*{problem*}{Problem}

\theoremstyle{fancy2}

\newtheorem{rem}[equation]{Remark}
\newtheorem*{rem*}{Remark}

\newcommand{\cref}[1]{Corollary~\ref{#1}}

\newcommand{\lref}[1]{Lemma~\ref{#1}}
\newcommand{\pref}[1]{Proposition~\ref{#1}}
\newcommand{\rref}[1]{Remark~\ref{#1}}

\newcommand{\e}{\epsilon}

\newcommand{\RP}{\mathbb{R\mkern1mu P}}
\newcommand{\CP}{\mathbb{C\mkern1mu P}}
\newcommand{\HP}{\mathbb{H\mkern1mu P}}
\newcommand{\CaP}{\mathrm{Ca}\mathbb{\mkern1mu P}^2}
\newcommand{\Sph}{\mathbb{S}}

\newcommand{\C}{{\mathbb{C}}}
\newcommand{\R}{{\mathbb{R}}}
\newcommand{\Z}{{\mathbb{Z}}}
\newcommand{\N}{{\mathbb{N}}}

\newcommand{\F}{\ensuremath{\operatorname{F}}}

\newcommand{\SO}{\ensuremath{\operatorname{SO}}}

\renewcommand{\O}{\ensuremath{\operatorname{O}}}
\newcommand{\Sp}{\ensuremath{\operatorname{Sp}}}
\newcommand{\U}{\ensuremath{\operatorname{U}}}
\newcommand{\SU}{\ensuremath{\operatorname{SU}}}
\newcommand{\Spin}{\ensuremath{\operatorname{Spin}}}

\newcommand{\T}{\ensuremath{\operatorname{T}}}
\renewcommand{\S}{\ensuremath{\operatorname{S}}}

\newcommand{\fg}{{\mathfrak{g}}}

\newcommand{\ft}{{\mathfrak{t}}}

\def\con#1=#2(#3){#1 \equiv #2 \bmod{#3}}

\newcommand{\tr}{\ensuremath{\operatorname{tr}}}
\newcommand{\diag}{\ensuremath{\operatorname{diag}}}

\DeclareMathOperator{\Id}{Id}

\newcommand{\no}{\noindent}
\newcommand{\co}{{cohomogeneity }}
\newcommand{\coo}{{cohomogeneity one }}

\newcommand{\ra}{\rangle}
\newcommand{\la}{\langle}

\newcommand{\brck}{{\rm [}}

\newcommand{\tz}{{\ft}_{\Z}}

\newcommand{\dd}{{\mathcal D}}

\newcommand{\spa}{\mbox{span}}
\newcommand{\kk}{\kappa}

\def\mod{{\rm\ mod \,}}

\def\qed{\ifhmode\unskip\nobreak\fi\ifmmode\ifinner\else\hskip5 pt
\fi\fi\hbox{\hskip5 pt \vrule width4 pt height6 pt depth1.5 pt
\hskip 1pt }}

\begin{document}

\title{Orbifold fibrations of Eschenburg spaces}  

\author{Luis A. Florit}
\address{IMPA: Est. Dona Castorina 110, 22460-320, Rio de Janeiro,
Brazil} \email{luis@impa.br}
\author{Wolfgang Ziller}
\address{University of Pennsylvania: Philadelphia, PA 19104, USA}
\email{wziller@math.upenn.edu}
\thanks{The first author was supported by CNPq-Brazil and the second author by
the Francis J. Carey Term Chair, the Clay Institute, by a grant from the
National Science Foundation, and by IMPA.}

\maketitle


Compact manifolds that admit a metric with positive sectional curvature
are still poorly understood. In particular, there are few known
obstructions for the existence of such metrics. By Bonnet-Meyers the
fundamental group must be finite, by Synge it has to be 0 or $\Z_2$ in
even dimensions, and the $\hat{A}$-genus must vanish when the manifold is
spin. For non-negative curvature, besides some results on the structure
of the fundamental group, we have Gromov's Betti number theorem which
states that they are bounded by a constant that only depends on the
dimension. In fact, there is no known obstruction that distinguishes the
class of simply connected manifolds which admit positive curvature from
the ones that admit non-negative curvature.

It is therefore surprising that there are very few known examples with
positive curvature. They all arise as quotients of a compact Lie group,
endowed with a left invariant metric, by a subgroup of isometries acting
freely. They consist, apart from the rank one symmetric spaces, of
certain homogeneous spaces in dimensions $6,7,12,13$ and $24$ due to
Berger \cite{Be}, Wallach \cite{Wa}, and Aloff-Wallach \cite{AW}, and of
biquotients in dimensions $6,7$ and $13$ due to Eschenburg
\cite{E1},\cite{E2} and Bazaikin \cite{Ba}.

A different method of searching for new positively curved examples is
suggested by another property that many (but not all) of the known
examples share: they are the total space of a fiber bundle where the
projection onto the base is a Riemannian submersion. It is therefore
suggestive to look for new examples which admit fiber bundle structures.
Weinstein (\cite{We}) studied this question by considering Riemannian
submersions with totally geodesic fibers, such that the sectional
curvatures spanned by a horizontal and a vertical vector are positive.
Even this weaker condition on a fiber bundle, called fatness, is already
strong (\cite{DR},\cite{Zi}).

The concept of fatness and any further curvature computations can be done
in a larger category of bundles, where the spaces involved are orbifolds
and the bundle structure is an orbifold one. Even if one is only
interested in manifolds this is an important generalization. In fact,
one can now give many of the other known examples of positive curvature
metrics an orbifold bundle structure as well; see Section 2 for a
summary. The only case where such an orbifold bundle structure was not
known, until now, is the family of (generic) Eschenburg biquotients
$E_{p,q}=\SU(3)/\!/\S^1$, where the $\S^1=\{z\in\C:|z|=1\}$ action on
$\SU(3)$ is given by
$$
z\cdot g=z^p\,g\,\overline z^{\,q},\;
$$
with $p,q\in\Z^3$, {\small$\sum$} $\!p_i=$ {\small$\sum$} $\!q_i$, and
$z^p:=\diag(z^{p_1}, z^{p_2}, z^{p_3})\in\U(3)$. This biquotient is an
orbifold if and only if $p-q_\sigma\neq 0$ for all permutations
$\sigma\in S_3$, where we have set
$q_\sigma=(q_{\sigma(1)},q_{\sigma(2)},q_{\sigma(3)})$.
If one has that $\gcd(p_i-q_{j}, p_{i'}-q_{j'})=1$, for all
$i\neq i',j\neq j'$, the quotient is a manifold. Further conditions must
be satisfied for an Eschenburg metric to have positive curvature; see
Section 1.

\bigskip

The purpose of this paper is to study all isometric circle actions on
Eschenburg manifolds (and more generally orbifolds) which act almost
freely, i.e., their isotropy groups are finite, and thus give rise to
principal orbifold bundle structures. One easily sees that they all
indeed admit such actions and we will examine in detail their geometric
properties. In particular, we obtain a large new family of 6-dimensional
orbifolds with positive sectional curvature and with small singular
locus.

\smallskip

By \cite{GSZ}, an isometric circle action on $E_{p,q}$ is given by a
biquotient action on $\SU(3)$ that commutes with the original one or,
equivalently, a $\T^2=\S^1\times \S^1$ biquotient action on $\SU(3)$ that
contains the original circle as a subgroup. So, given $a,b\in\Z^3$ with
\mbox{{\small$\sum$} $\!a_i=$ {\small$\sum$} $\!b_i$}, we define a circle
action $\S^1_{a,b}$ on $E_{p,q}$ by
$$
w\cdot [g]=[w^ag\,\overline w^{\,b}],\ \ \ w\in \S^1.
$$
The projection onto the quotient $\hat\pi:E_{p,q}\to O^{a,b}_{p,q}$
is then an orbifold principal bundle if this action
is almost free. We will show:
\begin{main}
The  circle action $\S^1_{a,b}$ on $E_{p,q}$ is almost free if and
only if
$$
(p-q_\sigma)\  \text{ and } \ (a-b_\sigma)\ \text{ are linearly
independent, for all }\, \sigma\in S_3 .
$$
The quotient $O^{a,b}_{p,q}$ is then an orbifold whose singular locus is
the union of at most nine orbifold 2-spheres and six points that are
arranged according to the schematic diagram in Figure 1.
\end{main}

\begin{picture}(100,160)  
\small \put(220,140){\line(-2,-1){50}}
\put(220,140){\line(2,-1){50}} \put(220,140){\line(0,-1){100}}
\put(220,40){\line(-2,1){50}} \put(220,40){\line(2,1){50}}
\put(270,115){\line(0,-1){50}} \put(170,115){\line(0,-1){50}}
\put(170,115){\line(2,-1){100}} \put(270,115){\line(-2,-1){100}}
\put(210,148){$C_{\Id}$} \put(148,54){$C_{(123)}$}
\put(148,122){$C_{(23)}$} \put(270,122){$C_{(12)}$}
\put(270,54){$C_{(132)}$} \put(210,25){$C_{(13)}$}
\put(240,40){${\mathcal L}_{13}$} \put(186,40){${\mathcal L}_{31}$}
\put(273,88){${\mathcal L}_{21}$} \put(150,88){${\mathcal L}_{23}$}
\put(184,135){${\mathcal L}_{11}$} \put(240,135){${\mathcal
L}_{33}$} \put(202,60){${\mathcal L}_{22}$} \put(188,107){${\mathcal
L}_{32}$} \put(244,95){${\mathcal L}_{12}$}
\put(217.3,137){$\bullet$} \put(217.3,37){$\bullet$}
\put(167.4,112){$\bullet$} \put(267.3,112){$\bullet$}
\put(267.3,62){$\bullet$} \put(167.4,62){$\bullet$} \put(108,0){\it
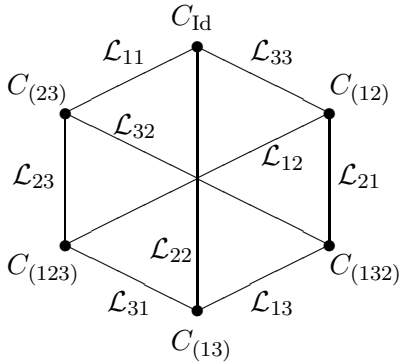
\small Figure 1. The structure of the singular locus.}
\end{picture}  

\bigskip

The lift of the orbifold 2-spheres to $\SU(3)$ consists of
the nine copies of $\U(2)$ inside $\SU(3)$ given by
$$
\U(2)_{ij}=\left\{ \tau_{i}\left[
\begin{array}{cc}
A & 0\\
0 & \det \overline A
\end{array}
\right]\tau_{j}\, : \ A \in \U(2) \right\}\; ,  \quad 1\leq
i,j\leq 3,
$$
where $\tau_r\in \O(3)$ is the linear map that interchanges the
$r^{\underline{th}}$ vector of the canonic basis with the third one.
The lift of the six singular points consists of the six copies of $\T^2$
given by
$$
\T^2_\sigma = \sigma \diag(z,w,\bar{z}\bar{w}),
$$
where $\sigma\in S_3\subset\O(3)$ is a permutation matrix. We define the
{\it parity} of each one of these six singular points to be the parity of
the corresponding $\sigma$. Clearly each $\U(2)$ contains two $\T^2$'s
and each $\T^2$ is contained in three $\U(2)$'s. They are also arranged
according to Figure 1, where edges correspond to the $\U(2)$'s and
vertices to the $\T^2$'s.

For the lift of the singular locus to $E_{p,q}$ under the fibration
$\hat{\pi}$, the edges in Figure~1 represent totally geodesic lens
spaces, and the vertices are closed geodesics. If $E_{p,q}$ is smooth,
the lens spaces are smooth as well. In Section 3 we will also determine
the isotropy groups corresponding to these lens spaces and closed
geodesics, which can also be interpreted as the orbifold groups of the
orbifold quotient. They are constant along each of these closed
geodesics, and along each lens space (outside the closed geodesics).

\medskip

In Section 4 we examine the question of how to minimize the singular
locus in $O^{a,b}_{p,q}$ and its orbifold groups. There exist some
Eschenburg spaces which admit a free circle action, but in general
the most one can hope for is an isometric circle action such that
the singular locus of the quotient consists of a single point. It
turns out that there is a topological obstruction to the existence
of such an action.

The most basic topological invariant of an Eschenburg space is the order
$h$ of the cyclic group $H^4(E_{p,q},\Z)$. We also associate to $E_{p,q}$
the integers (mod $h$) denoted by $\alpha(\sigma,\e_1,\e_2)\in\Z_h$,
where $\sigma\in S_3$ and $\e_1,\e_2=\pm1$, that only depend on $p,q$\,;
see \eqref{e:alpha} for an explicit formula. We show:

\begin{main}
Let $E_{p,q}$ be an Eschenburg manifold equipped with an Eschenburg
metric. Then, there exists an isometric circle action on $E_{p,q}$
whose singular locus is composed of at most 3 points with the same
parity if and only if $\alpha(\sigma,\e_1,\e_2)=0$ for some choice
of $\sigma\in S_3$, $\e_1,\e_2=\pm1$.
\end{main}

This implies in particular that a generic Eschenburg space does not admit
an isometric circle action with only one singular point. However, it is
easy to find examples for which this is the case; see Section 4. Observe
that all such examples have non-negative curvature, since this holds for
the Eschenburg metric in general.

On the other hand, for positive curvature the situation is different. To
illustrate this, we study in detail the case of general \coo Eschenburg
manifolds, that is, $E_d=E_{(1,1,d),(0,0,d+2)}$, $d\geq 0$, which have
positive curvature when $d>0$. For $d\leq 2$, it is known that $E_d$
admits a free isometric circle action, in fact even free actions by
$\SO(3)$ (\cite{Sh}). For the remaining cases, we prove the following.

\begin{main}
Let $E_d$ be any cohomogeneity one Eschenburg manifold, $d\geq 3$,
equipped with a positively curved Eschenburg metric. Then:
\begin{itemize}
\item[$i)$] There is no isometric $\ \S^1$ action on $E_d$ with only one
singular point.
\end{itemize}
In the following particular examples the singular locus of the
isometric circle action $\S^1_{a,b}$ on $E_d$ consists of:
\begin{itemize}
\item[$ii)$] Two points with equal orbifold groups $\Z_{d+1}$ if
$a=(0,-1,1)$ and $b=(0,0,0);$
\item[$iii)$] Two points with equal orbifold groups $\Z_{d-1}$ if
$a=(0,1,1)$ and $b=(2,0,0)$, and $3$ does not divide $d-1$.
If it does divide, we get in addition the orbifold 2-sphere
joining these two singular points, with orbifold group $\Z_3\,;$
\item[$iv)$]
A smooth totally geodesic 2-sphere with orbifold group $\Z_{d-1}$ if
$a=(0,1,1)$ and $b=(0,0,2)$.
\end{itemize}
\end{main}

Hence, an interesting open question is whether part $(i)$ in Theorem~C
holds more generally for all positively curved Eschenburg manifolds.
Computer experiments support an affirmative answer:
{\it There is no positively curved Eschenburg manifold with
$|H^4(E_{p,q},\Z)|\leq\ $100,000, altogether 103,569,197 spaces, that
admits an isometric circle action whose singular locus is a single point.}
If they do not exist, it would be interesting to understand the phenomena
behind this difference between non-negative and positive curvature.

Finally, observe that the most regular orbifold we obtain in
Theorem C is given by
$\diag(z,zw,z^3w)\backslash\SU(3)/
\diag(1,1,\overline z^{\,5}\overline w^{\,2})$,
which is a compact 6-dimensional positively curved orbifold which has
only two singular points with orbifold groups $\Z_2$.
\bigskip

We would like to thank C.\,G.\,Moreira for helpful discussions. This work
was done while the second author was visiting IMPA and he would like to
thank the Institute for its hospitality.

\vskip .8cm

\section{Preliminaries}  

Recall that an orbifold is a topological space which locally is the
quotient of an open set $U\subset\R^n$ under the effective action of a
finite group $\Gamma$ that fixes $p\in U$. The group $\Gamma$ is called
the orbifold group at the projection of $p$ in $U/\Gamma$, and the
required natural compatibility conditions for 2 overlapping orbifold
charts implies that the orbifold group is well defined. An orbifold
metric is a Riemannian metric on each chart $U$ such that $\Gamma$ acts
isometrically. In many ways orbifolds can be treated just like manifolds.
For the purpose of local geometric calculations, they can be done on the
smooth metric since all geometric objects are invariant under isometries.
The simplest examples of orbifolds are manifolds divided by a finite
group (so called good orbifolds). More generally, if a compact Lie group
$G$ acts isometrically on a Riemannian manifold $M$ such that all
isotropy groups are finite (so called almost free actions), then $M/G$ is
an orbifold, as follows immediately from the slice theorem for the group
action. Moreover, it also implies that the orbifold groups are the
isotropy groups, divided by the ineffective kernel. In our case,
orbifolds will be obtained as quotients of Eschenburg spaces under circle
actions with finite isotropy groups.

\medskip

We now introduce some notations that will be helpful and discuss general
properties of Eschenburg spaces. We denote the diagonal matrices
$\dd=\C^3\subset \C^{3\times 3}$ by
$x=(x_1,x_2,x_3)=\diag(x)\in \C^{3\times 3}$. For an element in the
symmetric group $S_3$ that takes $i\to j\to k \to i$, or $i\to j\to i$,
we use the notation $(i j k)$, or $(i j)$, respectively. We have a
natural action of $S_3$ on $\dd$ defined by
$x_\sigma=\diag(x_{\sigma(1)},x_{\sigma(2)},x_{\sigma(3)})$,
$\sigma\in S_3$, $x\in \dd$. If $z\in\S^1=\{z\in \C: |z|=1\}$ and
$p\in\Z^3\subset \R^3$, we denote by $z^p =
\diag(z^{p_1},z^{p_2},z^{p_3})\in \U(3)$. Observe that
$$
G=\{(g_1,g_2)\in \U(3)\times \U(3): \det g_1 = \det g_2\}
$$
acts on $\SU(3)$ by $(g_1,g_2)\cdot g = g_1\,g\,g_2^{-1}$. It has a
maximal torus $\T^5=(\dd\times\dd)\cap G$, whose Lie algebra is
$\ft=(\dd\times\dd)\cap \fg$. Let $\tz$ be the lattice in $\ft$ given by
$\tz={\ft}\cap (2\pi i\Z^3\times 2\pi i\Z^3)$, that we identify with
$$
\tz=\{(p,q)\in \Z^3\times \Z^3: \tr p =\tr q\}
$$
via $(p,q)\to (2\pi i p , 2\pi i q)$.

\bigskip

For each $(p,q)\in \tz$ we define an\ \ $\S^1$ action on $\SU(3)$ by
\begin{equation}\label{a1}
z\cdot g = z^p g\, \overline z^{\,q},\ \ \ \ z\in  \S^1,\ \ g\in \SU(3).
\end{equation}
This action is easily seen to be almost free if and only if
\begin{equation}\label{af}
p-q_\sigma\neq 0, \ \ \ \forall\ \sigma\in S_3,
\end{equation}
since this is clearly the case only when the action is free on the Lie
algebra level. In this situation, the quotient is a $7$-dimensional
orbifold, which we call the {\it Eschenburg orbifold} $E_{p,q}$, that
comes with the projection $\pi=\pi_{p,q}\colon \SU(3)\to E_{p,q}$\,.
Furthermore, the action is free if and only if
\begin{equation}\label{free}
\gcd(p_1-q_{\sigma(1)}, p_2-q_{\sigma(2)})=1,
\ \ \ \forall\ \sigma\in S_3,
\end{equation}
which is equivalent to $\gcd(p_i-q_{j}, p_{i'}-q_{j'})=1$, for all
$i\neq i',j\neq j'$. In this case, $E_{p,q}$ is a smooth
$7$--dimensional manifold, called the {\it Eschenburg manifold}
$E_{p,q}$. Finally, the action is effective only when
\begin{equation}\label{e}
\gcd(\{p_i-q_j:1\leq i,j \leq 3\})=1.
\end{equation}
Indeed, $z\in \S^1$ fixes $g\in\SU(3)$ if and only if $z^p$ is
conjugate to $z^q$ and if this is true for all~$g$, we necessarily
have $z^p=z^q=\lambda \Id$ for some $\lambda$. Again, one can take
$1\leq i \leq 2$ in the above since $\tr p =\tr q$.

\medskip

The {\it Eschenburg metric} on $E_{p,q}$ is the submersion metric
obtained by scaling the biinvariant metric on $\SU(3)$ in the direction
of $\U(2)_{jj}\subset\SU(3)$ for some $1\leq j \leq 3$, and it has
positive sectional curvature if and only if, for all $1\leq i \leq 3$,
\begin{equation}\label{pos}
p_i \notin [\min(q_1,q_2,q_3),\max(q_1,q_2,q_3)], \ \ \text{or}\ \
q_i \notin [\min(p_1,p_2,p_3),\max(p_1,p_2,p_3)].
\end{equation}
If this condition is satisfied, we will call $E_{p,q}$ a
{\it positively curved Eschenburg space}. Since the proof of this fact is
a Lie algebra computation, it still remains valid if we consider, more
generally, Eschenburg orbifolds. Thus in the orbifold category there
exists a much larger class of compact positively curved examples. In
fact, any six integers satisfying \eqref{pos} will determine one, since
\eqref{pos} implies \eqref{af}.

\medskip

Observe that, since $(p,q)$,$(-p,-q)$ and $(p+k\Id,q+k\Id)$, $k\in \Z$,
induce the same action and $\tr p= \tr q$, an Eschenburg orbifold is
determined by only 4 integers, that is, certain 4 of the differences
$p_i-q_j$. In fact,
$\{p_{\tau(1)}-q_{\sigma(2)},p_{\tau(1)}-q_{\sigma(3)},
p_{\tau(2)}-q_{\sigma(1)},p_{\tau(3)}-q_{\sigma(1)}\}$
defines an Eschenburg orbifold for fixed given permutations
$\tau,\sigma\in S_3$. Furthermore, notice that $E_{p,q}$ has positive
sectional curvature if and only if two rows or two columns of the matrix
$A_{ij}=p_i-q_j$ contain integers with the same sign, and it is a
manifold if and only if any two entries not in the same row or column are
relatively prime.

\medskip

There are two natural subclasses of Eschenburg manifolds. One is the
family of cohomogeneity two Eschenburg spaces corresponding to
$(p,q)=((c,d,e),(0,0,c+d+e))$ with $\gcd(c,d)=\gcd(d,e)=\gcd(e,c)=1$.
They admit an isometric action of $\T^2\times\SU(2)$ such that the
quotient is two-dimensional. A further subclass is the family of
cohomogeneity one Eschenburg spaces corresponding to
$(p,q)=((1,1,d),(0,0,d+2))$ which admit an isometric action of
$\SU(2)\times\SU(2)\times\S^1$ such that the quotient is one-dimensional.
In general, Eschenburg spaces admit an isometric action of $\T^3$ such
that the quotient is four-dimensional. In \cite{GSZ} it was shown that
these three groups are indeed the identity component of the isometry
group. In particular, the isometry has rank three in all cases.

\medskip

The only homological invariant that varies for different Eschenburg
manifolds is the order $h=h(E_{p,q})$ of the cohomology group
$H^4(E_{p,q},\Z)=\Z_h$. This integer is given by
$$
h=|p_1p_2+p_1p_3+p_2p_3-q_1q_2-q_1q_3-q_2q_3|
$$
(see \cite{E2}), which can be rewritten, up to sign, as
\begin{equation}\label{e:h}
h=(p_{\tau(1)}-q_{\sigma(2)})(p_{\tau(1)}-q_{\sigma(3)})-
(p_{\tau(2)}-q_{\sigma(1)})(p_{\tau(3)}-q_{\sigma(1)}),
\end{equation}
for any permutations $\tau,\sigma\in S_3$. Moreover, the integer $h$
must be odd for Eschenburg manifolds (see \cite{Kr}, Remark 1.4).
Notice also that, if $E_{p,q}$ is positively curved, we can assume
that $p_{\tau(1)}-q_{\sigma(2)}$, $p_{\tau(1)}-q_{\sigma(3)}$,
$p_{\tau(2)}-q_{\sigma(1)}>0$ and $p_{\tau(3)}-q_{\sigma(1)}<0$, and
hence there are only finitely many positively curved Eschenburg
manifolds for a given order $h$ (\cite{CEZ}).

\medskip

Finally, we introduce notations for a few orbifolds that we will need.
Given $p,q,d\in\Z$, $d\neq 0$, the {\it lens space} $L(p,q,d)$ is
the quotient
$$
L(p,q,d):=\Sph^3\!/\Z_d,
$$
where the action of $\Z_d=\{\xi\in \S^1:\xi^d=1\}$ on
$\Sph^3=\{(x,y)\in\C^2:|x|^2+|y|^2=1\}$ is given by
\begin{equation}\label{ac}
\xi\cdot (x,y)=(\xi^px,\xi^qy).
\end{equation}
This orbifold is a smooth manifold when $\gcd(p,d)=\gcd(q,d)=1$.
When there is no restriction on $\xi$, we get the
{\it weighted complex projective space}
$$
\CP^1[p,q]:=\Sph^3/\S^1,
$$
that is, the\, $\S^1$ action is still given by (\ref{ac}), for
$\xi\in\S^1$. For convenience, we will still call lens space the orbifold
$$
L(p,q,0):=\Sph^1\times\CP^1[p,q].
$$

\bigskip

\section{Known orbifold fibrations}  

\medskip

In this section we collect, for the convenience of the reader, the
known fibrations and orbifold fibrations where the total space is
one of the known compact simply connected positively curved manifolds and
the projection is a Riemannian submersion. Here we will leave out the
rank one symmetric spaces with their well known Hopf fibrations.

\bigskip

We start with the homogeneous examples $G/H$, which have been classified
in \cite{AW} and \cite{BB}. All except for one admit homogeneous
fibrations of the form $K/H\to G/H\to G/K$ coming from inclusions
$H\subset K\subset G$:

\bigskip
1. The Wallach flag manifolds each of which is the total space of the
following fibrations:

\bigskip

$\bullet$ \centerline{$\Sph^2\to \SU(3)/\T^2 \to \CP^2$,}

\bigskip

$\bullet$ \centerline{$\Sph^4\to \Sp(3)/\Sp(1)^3 \to \HP^2$,}

\bigskip

$\bullet$ \centerline{$\Sph^8\to \F_4/\Spin(9) \to \CaP$.}

\bigskip

2. The Aloff-Wallach examples
$E_{0,q}=W_{q_1,q_2}=\SU(3)/\diag(z^{q_1},z^{q_2},\bar{z}^{q_1+q_2})$,
which have positive curvature when $q_1q_2(q_1+q_2)\ne 0$, admit
two kinds of fibrations:

\bigskip
$\bullet$ \centerline{$\Sph^1\to W_{q_1,q_2}\to \SU(3)/\T^2,$}
\bigskip
and a lens space fibration

$\bullet$ \centerline{$\Sph^3/\Z_{q_1+q_2}\to W_{q_1,q_2} \to \CP^2,$}
\smallskip

\no where the fiber is
$\U(2)/\diag(z^{q_1},z^{q_2})=\SU(2)/\Z_{q_1+q_2}$.

\smallskip

3. The Berger example $\SU(5)/\Sp(2)\cdot \S^1$ admits a
fibration

\bigskip
$\bullet$ \centerline{$\RP^5\to\SU(5)/\Sp(2)\cdot \S^1\to\CP^4,$}
\smallskip

\no where the fiber is $\U(4)/\Sp(2)\cdot \S^1=\SU(4)/\Sp(2)\cdot
\Z_2 =\SO(6)/\O(5)=\RP^5$.

\smallskip

4. Finally for the homogeneous category, we have the Berger space
$\SO(5)/\SO(3)$. This space is special since $\SO(3)$ is maximal in
$\SO(5)$ and hence does not admit a homogeneous fibration. In \cite{GKS}
it was shown that $\SO(5)/\SO(3)$ is diffeomorphic to the total space of
an $\Sph^3$ bundle over $\Sph^4$, but the fibration is not a Riemannian
submersion of the positively curved metric. It was observed though by K.
Grove and the last author that the subgroup
$\SU(2)\subset\SO(4)\subset\SO(5)$ acts with only finite isotropy
groups and hence gives rise to an orbifold fibration

\smallskip
$\bullet$ \centerline{$\Sph^3\to\SO(5)/\SO(3)\to \Sph^4.$}
\smallskip

\no To see this, one observes that the action by $\SO(4)\subset
\SO(5)$ has \coo and from the group diagram of this action (see
\cite{GWZ}) it follows that there is a codimension two submanifold,
one of the singular orbits $\SO(4)/\O(2)$, along which the isotropy
group is a cyclic group $\Z_3$ and away from this orbit, the action
is free. The quotient is homeomorphic to $\Sph^4$ and the metric is
smooth, except along a Veronese embedding $\RP^2\subset\Sph^4$ where
the metric has an angle $2\pi/3$ normal to $\RP^2$.

\smallskip

The remaining known examples of positively curved compact manifolds
are given by biquotients:
\smallskip

5. For the Eschenburg spaces, there are 3 subfamilies which are known to
admit fibrations. One is the family of Aloff-Wallach spaces discussed
above. A second one arises from a free action of $\T^2$ whose quotient
is the inhomogeneous positively curved flag manifold discovered by
Eschenburg,
$\SU(3)/\!/\T^2:=\diag(z,w,zw)\backslash\SU(3)/\diag(1,1,z^2w^2)^{-1}$,
which gives rise to a fibration

\smallskip
$\bullet$ \centerline{$\Sph^1\to E_{p,q}\to \SU(3)/\!/\T^2,$}
\smallskip

\no for every $(p,q)$ of the form
$(p,q)=((p_1,p_2,p_1+p_2),(0,0,2p_1+2p_2))$.

\smallskip

The third subfamily consists of the \co two Eschenburg manifolds
defined by $(p,q)=((p_1,p_2,p_3),(0,0,\bar{p})$ where
$\bar{p}=p_1+p_2+p_3$ and the $p_i$'s are pairwise relatively prime.
They admit an action by $\SU(2)=\{\diag(A,1): A\in \SU(2)\}$
acting on the right since it commutes with the circle action. As was
observed in \cite{BGM}, it gives rise to an orbifold fibration

\bigskip
$\bullet$ \centerline{$F\to E_{p,q}\to
\CP^2[p_2+p_3,p_1+p_3,p_1+p_2]\ ,$}
\smallskip

\no where the fiber $F$ is $\RP^3$ if all $p_i$'s are odd, and $F=\Sph^3$
otherwise. Here the base is a \mbox{2-dimensional} weighted complex
projective space. Indeed, if one first uses the identification
$\SU(3)/\SU(2)\cong \Sph^5$ given by $[g]\mapsto g(e_3)$, the remaining
circle action becomes
$(v_1,v_2,v_3)\to (z^{\bar{p}-p_1}v_1,
z^{\bar{p}-p_2}v_2,z^{\bar{p}-p_3}v_3)$.
This also shows that the coordinate points in the weighted projective
space correspond to $\U(2)$'s inside $\SU(3)$, and the isotropy groups
are cyclic of order $\bar{p}-p_i$. On the other hand, we also have that
$\gcd(\bar{p}-p_i,\bar{p}-p_j)=:a>1$ for at least one pair $i,j$ and
hence the orbifold set also contains at least one $\CP^1\subset\CP^2$
with orbifold group $\Z_a$.

\smallskip

6. Finally, we have the Bazaikin biquotients
$B_q=\SU(5)/\!/\Sp(2)\cdot \S^1$ given by
$$
B_q = \diag (z^{q_1} , \dots , z^{q_5} ) \backslash \SU(5)/ \diag
(A,z^{\,q_0})^{-1},
$$
where $A\in \Sp(2)\subset \SU(4)\subset\SU(5)$,
$q=(q_1,\cdots ,q_5)$ is an ordered set of odd integers and
${q_0}=\sum q_i$. In this case we obtain the fibration

\bigskip

$\bullet$ \centerline{$\RP^5\to B_q\to
\CP^4[q_0-q_1,\cdots,q_0-q_5].$}
\smallskip

\no by enlarging the biquotient action of $\Sp(2)\cdot \S^1$ to one
of $\SU(4)\cdot \S^1$. The isotropy groups and the weights are
obtained as in the case 5 above.

\medskip

It is an interesting fact that there are no other two tori that act
freely on $\SU(3)$ besides the ones in fibrations 2 and 5; see \cite{E2}.
It thus only remains the general family of Eschenburg spaces, which was
not known, until now, to admit an orbifold fibration.

\medskip
\section{Circle orbifold fibrations}  
\medskip

We now search for an almost free $\ \S^1$ action on $E_{p,q}$. We also
require that the circle action acts isometrically in the positively
curved Eschenburg metrics. As mentioned in Section 1, the isometry group
of a positively curved Eschenburg space has rank 3 and hence any circle
action is conjugate to one lying in a maximal torus. This maximal torus
can be chosen to be the 3-torus induced by the biquotient action of the
maximal torus $\T^5\subset G\subset \U(3)\times \U(3)$. We are thus
forced to consider circle actions induced by a circle inside this
3-torus. This amounts to finding an $\S^1$ action on $\SU(3)$ that
commutes with the one that defines $E_{p,q}$, in such a way that they
give together a $\T^2= \S^1\times \S^1\subset \T^5$ almost free action on
$\SU(3)$.

To describe this $\T^2\subset \T^5$, let $(a,b)\in \tz$, and write
the $\T^2$ action as
\begin{equation}\label{t2}
(z,w)\cdot g = w^az^p g\, \overline z^{\,q}\overline w^{\,b},\ \ \ \
z,w\in\S^1,\ g\in \SU(3).
\end{equation}
The action is almost free if and only if it is free at the Lie algebra
level. Since the Lie algebra of $\T^2$ is spanned by $i(p,q)$ and
$i(a,b)\in \mathfrak{u}(3)\times\mathfrak{u}(3)$, this holds if and only
if there are no $x,y\in\R$ such that $xp+ya$ is conjugate to $xq+yb$.
Since both are diagonal matrices, the almost free property is then
equivalent to
\begin{equation}\label{aff}
(a-b_\sigma)\ \ {\rm and}\ \ (p-q_\sigma)\ \
{\rm are\ \ linearly\ \ independent},\ \forall\ \sigma \in S_3.
\end{equation}
This $\T^2$ action on $\SU(3)$ defines a circle action on $E_{p,q}$,
which we denote by $\S^1_{a,b}$ and its quotient by $O_{p,q}^{a,b}$. This
circle action is clearly almost free if and only if the $\T^2$ action is
almost free, and, in this case, $O_{p,q}^{a,b}$ is an orbifold. Recall
that we have the projection ${\pi}_{p,q}\colon\SU(3)\to E_{p,q}$ and we
define a further projection
$\hat{\pi}_{p,q}^{a,b}\colon E_{p,q}\to O_{p,q}^{a,b}$.
If clear from context, we also denote these projections simply by $\pi$
and $\hat{\pi}$, respectively. Our purpose is to study the geometry of
the orbifold fibration
\begin{equation*}
\S^1\to E_{p,q}\to O_{p,q}^{a,b}.
\end{equation*}

 From now on we assume that $(p,q)$ and $(a,b)$ satisfy (\ref{aff}).
Notice that this implies that (\ref{af}) holds for $n(p,q)+m(a,b)$, for
all $(n,m)\in\Z^2\setminus \{0\}$. Thus $E_{a,b}$ is also an Eschenburg
orbifold and we obtain a symmetry in the process and the commutative
diagram of orbifold fibrations given in Figure $2$.

\begin{picture}(100,170)  
\small
\put(206,98){$\SU(3)$}
\put(210,90){\vector(-1,-1){30}}
\put(230,90){\vector(1,-1){30}}
\put(256,46){$E_{a,b}$}
\put(162,46){$E_{p,q}$}
\put(189.5,0){$O_{p,q}^{a,b} \ = \ O_{a,b}^{p,q}$}
\put(260,37){\vector(-2,-3){15}}
\put(180,37){\vector(2,-3){15}}
\put(296,98){$\S^1$}
\put(136,98){$\S^1$}
\put(295.2,90){\vector(-2,-3){20}}
\put(144.7,90){\vector(2,-3){20}}
\put(256,140){$\S^1$}
\put(173,140){$\S^1$}
\put(254.9,135){\vector(-1,-1){22}}
\put(185,135){\vector(1,-1){22}}
\put(215,140){$\T^2$}
\put(195,140){$\subset$}
\put(235,140){$\supset$}
\put(220,135){\vector(0,-1){20}}
\put(220,88){\vector(0,-1){70}}
\put(100,-25){\it Figure 2. Orbifold fibrations of Eschenburg spaces.}
\end{picture}  

\vspace{40pt}

Since $\S^1_{a,b}$ is almost free, there are only regular orbits
and exceptional orbits. We define the {\it exceptional set}\ \
${\mathcal S}^{a,b}_{p,q}\subset E_{p,q}$\ \ to be the union of all
exceptional orbits, which thus coincides with the set of points in
$E_{p,q}$ where the action is not free.

\medskip

Recall that we have the nine embeddings $\U(2)_{ij}\subset\SU(3)$,
$1\leq i,j\leq 3$, and the six two-dimensional tori
$\T^2_\sigma$, $\sigma\in S_3$. They give rise to their respective
images in $E_{p,q}$,
$$
{\mathcal L}_{ij}:=\pi(\U(2)_{ij}),\ \ \ \
\ C_\sigma:=\pi(T^2_\sigma).
$$
While the $C_\sigma$'s are clearly circles, we claim that the
${\mathcal L}_{ij}$'s are lens spaces. Indeed, if $g\in \U(2)_{ij}$
and $z\in S^1$, we get
\begin{equation}\label{lens}
\tau_iz^pg\overline z^{\,q}\tau_j= \left[
\begin{array}{ccc}
z^{p_{i_1}-q_{j_1}}x & z^{p_{i_1}-q_{j_2}}y & 0\\
-z^{p_{i_2}-q_{j_1}}\lambda\overline y &
z^{p_{i_2}-q_{j_2}}\lambda\overline x & 0\\
0 & 0 & z^{p_i-q_j}\overline{\lambda}
\end{array}
\right],
\end{equation} where $\lambda\in S^1$ and $(x,y)\in \Sph^3$ and we have
used the index convention $\{i_1,i_2,i\}=\{j_1,j_2,j\}=\{1,2,3\}$.
Taking a representative $g\in \SU(2)=\Sph^3$ in the orbit (i.e.
$\lambda=1$) and identifying the upper $2\times 2$ matrix in
(\ref{lens}) with its first row, we conclude that
${\mathcal L}_{ij}$ is the lens space
$$
{\mathcal L}_{ij}=L(p_{i_1}-q_{j_1},p_{i_1}-q_{j_2},p_i-q_j).
$$
 From their very definition, we see that each lens space
${\mathcal L}_{ij}$ contains precisely two of the circles
$C_\sigma\subset{\mathcal L}_{ij}$, where $\sigma$ is one of the two
permutations that satisfy $\sigma(i)=j$. Moreover, each circle is then
obtained as the intersection of three lens spaces. They are arranged as
shown in Figure $1$ in the Introduction.

\bigskip

We now proceed to investigate the structure of
${\mathcal S}_{p,q}^{a,b}$. The isotropy group $ \S^1_{[g]}\subset \S^1$
at $[g]\in E_{p,q}$ is the finite cyclic group given by the elements
$w\in\S^1$ such that there is $z\in \S^1$ with
\begin{equation}\label{sing}
g^{-1}z^pw^ag = z^qw^b.
\end{equation}
If $1\neq w \in  \S^1_{[g]}$, there must be $z\in \S^1$ and
$\sigma \in S_3$ such that
\begin{equation}\label{premain}
z^pw^a =z^{q_\sigma}w^{b_\sigma}.
\end{equation}
If there exists a $w$ such that the two matrices in \eqref{premain} have
a triple eigenvalue, then it acts trivially on all of $E_{p,q}$ and hence
belongs to the ineffective kernel, whose order we denote by $\kk_0\in\N$.

Otherwise, there are two possibilities. If all eigenvalues are distinct,
\eqref{sing} implies that $g(e_i) = \mu_ie_{\sigma^{-1}(i)}$ with
$\mu_1\mu_2\mu_3=1$ and thus $[g]\in C_\sigma$. Furthermore, it follows
that there exists a $\kk_\sigma\in\N$ such that
$S^1_{[g]}=\Z_{\kk_\sigma}$ for all $[g]\in C_\sigma$, i.e. the isotropy
group has constant order along the exceptional subset $C_\sigma$.

If, on the other hand, the matrices in (\ref{premain}) have a double
eigenvalue, the third eigenvalue must coincide as well since $(p,q)$ and
$(a,b)$ belong to $\tz$. Thus there exist $i,j$ with $1\leq i,j\leq 3$
such that $\overline z^{\,p_{i'}-q_{j'}}=w^{a_{i'}-b_{j'}}$, for all
$i'\neq i$, $j'\neq j$. It follows that $[g]\in {\mathcal L}_{ij}$, and
that $w$ fixes all of ${\mathcal L}_{ij}$. Hence there exists a
$\kk_{ij}\in\N$ such that the isotropy group along ${\mathcal L}_{ij}$ is
$\Z_{\kk_{ij}}$, except along the two circles
$C_\sigma\subset{\mathcal L}_{ij}$ where it is $\Z_{\kk_\sigma}$, with
$\sigma$ being one of the two permutations that satisfy $\sigma(i)=j$.

We thus have shown that the exceptional set is given by

\begin{equation}\label{e:struc}
{\mathcal S}^{a,b}_{p,q}\ \ =\bigcup_{(i,j)\in \Gamma}{\mathcal L}_{ij}\
\bigcup_{\sigma\in\Sigma}C_\sigma\ ,
\end{equation}
where $\Gamma=\{(i,j): \kk_{ij}>\kk_0,\ 1\leq i,j\leq 3\}$ and
$\Sigma=\{\sigma\in S_3:\kk_\sigma>\kk_0\}$.

\medskip

Notice that each lens space ${\mathcal L}_{ij}$ in this exceptional set
must be totally geodesic. To see this, it is sufficient to show that
$\U(2)_{ij}$ is totally geodesic in $\SU(3)$ with respect to the
Eschenburg metric since $\U(2)_{ij}$ is invariant under the circle action
$\S^1_{a,b}$. This in turn follows since $\U(2)_{ij}$ is the fixed point
set of an isometry on $\SU(3)$ of the form $g\to z^vg\bar{z}^{v_\sigma}$
for some $z^v$ with two equal diagonal entries and some permutation
$\sigma$. The circles $C_{\sigma}$, being intersections of two totally
geodesic submanifolds, are thus closed geodesics.

The metric on the lens spaces are in general not homogeneous, in fact
they are homogeneous if and only if the circle action on $\U(2)_{ij}$ is
one sided; see the proof of Theorem 4.1 in \cite{GSZ}. This is only
possible when the Eschenburg space has \co two. In the case where the
lens space is homogeneous, the metric is a Berger type metric, i.e.
induced from a metric on $\Sph^3$ shrunk in direction of the Hopf fibers.
In the \co one case there are 6 such homogeneous lens spaces and in the
remaining \co two spaces 3 of them are homogeneous. When the lens space
is not homogeneous, its isometry group still contains a 2-torus.

\bigskip

We next determine the order of the isotropy groups of our circle
actions. To do so, for $v\in \Z^n$, set
$\gcd(v)=\gcd(\{v_1,\dots,v_n\})$ and define
$$
\kappa(v,w):=\gcd(v)^{-1}\gcd(\{v_iw_j-v_jw_i: 1\leq i<j\leq n\}).
$$

\begin{prop}\label{orders}
For the almost free circle action $\,\S^1_{a,b}$ on $E_{p,q}$ given by
\eqref{t2} it holds that:
\begin{itemize}
\item[$(a)$] The ineffective kernel is $\Z_{\kk_0}$, $\kk_0=\kappa(P,A)$,
where $P,A\in\Z^6$ are the vectors whose components are $p_i-q_j$ and
$a_i-b_j$, respectively, with $i\neq j$ (the same index ordering for
both).
\item[$(b)$] The isotropy group along $C_\sigma$ is $\Z_{\kk_\sigma}$,
where
$$
\kk_\sigma=\frac{|(p_1-q_{\sigma(1)})(a_2-b_{\sigma(2)})-
(p_2-q_{\sigma(2)})(a_1-b_{\sigma(1)})|}
{\gcd(p_1-q_{\sigma(1)},p_2-q_{\sigma(2)})}.
$$
\item[$(c)$] The isotropy group along ${\mathcal L}_{ij}$, outside
$C_\sigma\subset{\mathcal L}_{ij} $ with $\sigma(i)=j$, is
$\Z_{\kk_{ij}}$, $\kk_{ij}=\kappa(V,W)$, where $V,W\in\Z^4$ are the
vectors whose components are $p_{i'}-q_{j'}$ and $a_{i'}-b_{j'}$,
respectively, with $i'\ne i , j'\ne j$ (the same index ordering for
both).
\end{itemize}
\end{prop}

\begin{proof}
We need the next elementary lemma concerning the lattice points
$\Z^n$ inside the parallelogram spanned by $v,w\in\Z^n$,
${\mathcal P}_{v,w}=\{t v + s w: t,s \in \brck 0,1)\}$.

\begin{lem}\label{parallelogram}
If $p=(p_1,\dots,p_n), a=(a_1,\dots,a_n)\in\Z^n$ are linearly
independent, the projection of the lattice points
${\mathcal P}_{p,a}\cap\Z^n$ inside ${\mathcal P}_{p,a}$ to the
$s$-coordinate is the set
$\Z_\kappa=\{i/\kappa: i=0,\dots,\kappa-1\}\subset [0,1)$, with
$\kappa=\kappa(p,a)$.
\end{lem}
\begin{proof}
We start with the case of $n=2$. The number of lattice points $\Z^2$
inside ${\mathcal P}_{p,a}$ is equal to its area, that is,
$\# ({\mathcal P}_{p,a}\cap\Z^2)=|p_1a_2-p_2a_1|$.
Indeed, using translations, we can assume that $p_i,a_i\geq0$ and,
if we inscribe ${\mathcal P}_{p,a}$ inside the rectangle
${\mathcal P}_{(p_1+a_1,0),(0,p_2+a_2)}$, it follows that
$\pm\# ({\mathcal P}_{p,a}\cap\Z^2)=
(p_1+a_1)(p_2+a_2)-a_1a_2-p_1p_2-2p_2a_1=p_1a_2-p_2a_1$.

Setting $p'=d^{-1}p$, observe that ${\mathcal P}_{p,a}$ is the union
of the $d=\gcd(p_1,p_2)$ disjoint rectangles
$kp' + {\mathcal P}_{p',a}$, $k=0,\dots,d-1$, since
$tp+s_0a,t'p+s_0a\in\Z^2$ implies that $(t-t')p\in \Z$. Thus the number
of points in the projection of ${\mathcal P}_{p,a}\cap\Z^2$ to the
$s$-coordinate is ${\gcd(p_1,p_2)^{-1}}{|p_1a_2-p_2a_1|}$. This proves
our claim if $n=2$.

For $n>2$, by the above argument we can assume $\gcd(p)=1$, and thus the
number of points of the projection of the lattice to the $s$-coordinate
is equal to $\kappa= \#({\mathcal P}_{p,a}\cap\Z^n)$. In the plane
$\Pi=\spa\{p,a\}$, consider the lattice $L=\Pi\cap\Z^n$. From the case
$n=2$ it follows that $\kk=|p\wedge a|/|v\wedge w|$, where $\{v,w\}$ is a
base of $L$, or, equivalently,
$\kk=\max \{|p\wedge a|/|v\wedge w|:
v,w\in L \text{\ are linearly independent}\}$.
Since $\gcd(p)=1$, we can take $v=p$ and thus $a=rp\pm\kappa w$, for
some $r\in\Z$. This implies that $\kk$ divides $(a_i-rp_i)p_j$ and
$(a_j-rp_j)p_i$ and thus $a_ip_j-a_jp_i$ for every $i,j$, i.e., $\kk$
divides $\kk(p,a)$. On the other hand, if we choose $u\in\Z^n$ with
$\la u,p\ra=1$, $\kk(p,a)$ divides
$\sum_ju_j(a_ip_j-a_jp_i)=a_i-\la u,a\ra p_i$ for every $i$, and hence
$w'=\kk(p,a)^{-1}(a-\la u,a\ra p)\in L$. Therefore,
$\kk\geq |p\wedge a|/|p\wedge w'|=\kk(p,a)$, and the lemma follows.
\end{proof}

Assume the action is not effective. Then, there exists $1\neq w\in \S^1$
such that, for every $g\in \SU(3)$, there exists $z=z_g\in \S^1$ with
$g^{-1}z_g^pw^ag=z_g^qw^b$. By choosing different $g$, it is easy to
see that there exist $z,\lambda\in \S^1$ such that
$z^pw^a=z^qw^b=\lambda \Id$. We can write this as
$\overline z^{\,(p_i-q_j)}=w^{(a_i-b_j)}$, for all $1\leq i,j\leq 3$.
If we set
\begin{equation}\label{zw}
z=e^{2\pi i t},\ \ w=e^{2\pi i s}, \ \ (t,s)\in [0,1)\times [0,1),
\end{equation}
this means that $tP+sA \in \Z^6$ and the claim follows from
\lref{parallelogram}.

 From our discussion above it follows that the determination of the
isotropy groups falls into 2 cases, depending on whether the matrices in
\eqref{premain} have simple eigenvalues or a double eigenvalue.

{\it Case 1: Simple eigenvalues.}
Using (\ref{zw}), there is $t\in[0,1)$ such that
\begin{equation}\label{main}
t(p-q_\sigma)+s(a-b_\sigma) \in \Z^3.
\end{equation}
We need only to consider the first two coordinates since
$(p,q),(a,b)\in \tz$ and hence we look for the lattice points $\Z^2$
inside the two dimensional parallelogram determined by
$v_\sigma=(p_1-q_{\sigma(1)},p_2-q_{\sigma(2)}),
w_\sigma=(a_1-b_{\sigma(1)},a_2-b_{\sigma(2)})$ and its projections to
the $s$-coordinate. Thus the claim follows from
\lref{parallelogram}.

{\it Case 2: Double eigenvalue.}
Since $(p,q)$ and $(a,b)$ belong to $\tz$, the third eigenvalue must
coincide as well. Thus there are \mbox{$1\leq i,j\leq 3$} such that
$\overline z^{\,p_{i'}-q_{j'}}=w^{a_{i'}-b_{j'}}$, for all $i'\neq i$,
$j'\neq j$. This is equivalent to $tV+sW \in \Z^4$ and we apply
\lref{parallelogram}.
\end{proof}

\begin{rem}\label{r:afcond}
Notice that, by \pref{orders} ($b$) and \eqref{aff}, the circle action
given by \eqref{t2} is almost free if and only if $\kk_\sigma\neq 0$ for
all $\sigma\in S_3$. Observe also that the order of the isotropy groups
are in terms of the possibly ineffective action. The orders need to be
divided by $\kk_0$ to obtain the isotropy groups of the action when made
effective. In explicit computations it is useful to notice that since
$(p,q),(a,b)\in \tz$, the number of entries in the definition of
$\kappa(P,A)$ in part $(a)$ can be reduced from 6 to 4 numbers, and for
$\kk_{ij}$ from 4 to 3.
\end{rem}

We finally discuss the singular locus of the orbifold $O_{p,q}^{a,b}$,
i.e., $\hat{\pi}({\mathcal S}^{a,b}_{p,q})\subset O_{p,q}^{a,b}$.
Each $\hat{\pi}(C_\sigma)$ will be called a {\it vertex} of
$\hat{\pi}({\mathcal S}^{a,b}_{p,q})$, while
$\hat{\pi}({\mathcal L}_{ij})$ will be called a {\it face}.

For this purpose, we assume for simplicity that the Eschenburg space
$E_{p,q}$ is a smooth manifold. In this case (\ref{free}) implies that
the lens spaces ${\mathcal L}_{ij}$ are all smooth manifolds as well.
This is evident if all $p_i$ are distinct from $p_j$. If two of them
agree, ${\mathcal L}_{ij}$ is either equal to $L(0,a,1)=\Sph^3$ or to
$L(1,1,0)=\Sph^2\times\Sph^1$ (the latter two are not possible in
positive curvature).

Clearly, by \pref{orders} and the slice theorem, each
$\hat{\pi}(C_\sigma)$ is a point with orbifold group
$\Z_{\kk_\sigma/\kk_0}$, while the orbifold group of
$\hat{\pi}({\mathcal L}_{ij})$ is $\Z_{\kk_{ij}/\kk_0}$,
outside $\hat{\pi}(C_\sigma)$ and $\hat{\pi}(C_{\sigma'})$
with $\sigma(i)=\sigma'(i)=j$.

On the other hand, each face $\hat{\pi}({\mathcal L}_{ij})$ is two
dimensional and is itself an orbifold which is totally geodesic in
$O_{p,q}^{a,b}$. We claim that it is homeomorphic to $\Sph^2$ and that it
has only two orbifold points, namely the vertices
$\hat{\pi}(C_\sigma),\hat{\pi}(C_{\sigma'})
\in\hat{\pi}({\mathcal L}_{ij})$
with orbifold angles $2\pi \kk_{ij}/\kk_\sigma$ and
$2\pi\kk_{ij}/\kk_{\sigma'}$, respectively. Indeed, $\S^1_{a,b}$
preserves ${\mathcal L}_{ij}$, with $C_\sigma$ and $C_{\sigma'}$ as two
of its orbits. We now apply the slice theorem of the action
restricted to ${\mathcal L}_{ij}$ at a point in $C_\sigma$ where the
isotropy group is $\Z_{\kk_\sigma}$. This isotropy group acting on the
two dimensional slice has $\Z_{\kk_{ij}}$ as its ineffective kernel,
while the quotient group $\Z_{\kk_\sigma}/\Z_{\kk_{ij}}$ acts effectively
via a finite rotation group. Thus, $\Z_{\kk_\sigma/\kk_{ij}}$ is the
orbifold group of $\hat{\pi}(C_\sigma)$ as a singular point of the two
dimensional orbifold $\hat{\pi}({\mathcal L}_{ij})$. Away from
$C_\sigma$, $C_{\sigma'}\subset {\mathcal L}_{ij}$, the circle action is
free, modulo its ineffective kernel $\Z_{\kk_{ij}}$, and hence
$\hat{\pi}({\mathcal L}_{ij})$ has only two orbifold points. As a
consequence, the singular two sphere $\hat{\pi}({\mathcal L}_{ij})$ is
smooth if and only if $\kk_{ij}=\kk_{\sigma}=\kk_{\sigma'}$. Since the
metric on the lens space has $\T^2$ as its isometry group, the orbifold
still admits a circle of isometries, i.e., it is rotationally symmetric.

\section{Minimizing the singular set}  

Given an Eschenburg space $E_{p,q}$, it is natural to try to find
quotients $O=O_{p,q}^{a,b}$ as regular as possible by minimizing the
exceptional set ${\mathcal S} ={\mathcal S}^{a,b}_{p,q}$. This can be
achieved in several ways, e.g., by making the orders $\kk_\sigma$ and
$\kk_{ij}$ as small as possible, or $\kk_{ij}=1$ for most pairs
$1\leq i,j\leq 3$ to eliminate the corresponding lens spaces from
${\mathcal S}$.

\medskip

Recall that there exist precisely two 2-tori
$\T^2\subset \U(3)\times \U(3)$ that act freely on $\SU(3)$ as a
biquotient action. In particular, there are two infinite
families of Eschenburg spaces which admit a free $\S^1$ action; see
Section 2. For the general case, we can try to minimize the singular set
as follows.

\begin{prop}\label{p:3lens}
For each Eschenburg manifold $E_{p,q}$ endowed with an Eschenburg metric,
there exists an isometric circle action whose exceptional set is composed
of at most 3 totally geodesic lens spaces, intersecting along one closed
geodesic, and the order of the cyclic isotropy groups of these lens
spaces is bounded by $h=|H^4(E_{p,q},\Z)|$.
\end{prop}

\begin{proof}
Fix $\sigma \in S_3$, and let $\sigma'=\sigma\circ(123)$,
$\sigma''=\sigma\circ(132)\in S_3$ be the two permutations with the
same parity of $\sigma$. Observe that if $\kk_{\sigma'}=1$, by
\eqref{e:struc} and \pref{orders} the three lens spaces that contain
$C_{\sigma'}$ consist of regular points. If, in addition,
$\kk_{\sigma''}=1$, then 6 of the lens spaces are regular, and therefore
${\mathcal S}$ would be composed of at most the 3 lens spaces that meet
at the circle $C_{\sigma}$. In particular, in this situation all the
$\kk_{ij}$ must divide $\kk_\sigma$. Now, we claim that this can always
be done, and with $\kk_\sigma\leq h$.

Fix $\e_1,\e_2=\pm 1$. Since $E_{p,q}$ is smooth, there are
$x,y,z,w\in\Z$ such that
\begin{equation}\label{e:cofactors}
\left\{\begin{array}{l}
x(p_1-q_{\sigma(2)})-y(p_2-q_{\sigma(3)})=\e_1,
\vspace{2ex}\\
w(p_1-q_{\sigma(3)})-z(p_2-q_{\sigma(1)})=\e_2.
\end{array} \right.
\end{equation}
In view of (\ref{free}), the set of all solutions of
\eqref{e:cofactors} is given by $x'=x+k_1(p_2-q_{\sigma(3)})$,
$y'=y+k_1(p_1-q_{\sigma(2)})$, $w'=w+k_2(p_2-q_{\sigma(1)})$ and
$z'=z+k_2(p_1-q_{\sigma(3)})$, with $k_1,k_2\in\Z$. Now, we define
$(a,b)\in \tz$, which depends on $\sigma, \e_1,\e_2,k_1,k_2$, by
$a=(-z',-x',y'+w')$, $b_{\sigma}=(w'-x',y'-z',0)$. By definition
this implies that $\kk_{\sigma'}=\kk_{\sigma''}=1$ for
$\S^1_{a,b}$. The orders of the other 4 circles are given by
\begin{equation}\label{e:orders}
\left\{\begin{array}{l}
\ \ \kk_{\sigma}\ \ \ =|s h +\ (x+y-z)(p_1-q_{\sigma(1)})
-(w+z-x)(p_2-q_{\sigma(2)})|,
\vspace{1ex}\\
\kk_{\sigma\circ(12)}=|s(p_1-q_{\sigma(2)})(p_2-q_{\sigma(1)})-
w(p_1-q_{\sigma(2)})+y(p_2-q_{\sigma(1)})|,
\vspace{1ex}\\
\kk_{\sigma\circ(23)}=|s(p_2-q_{\sigma(3)})(p_3-q_{\sigma(2)})+
(z+w)(p_2-q_{\sigma(3)})+x(p_3-q_{\sigma(2)})|,
\vspace{1ex}\\
\kk_{\sigma\circ(13)}=|s(p_1-q_{\sigma(3)})(p_3-q_{\sigma(1)})-
(x+y)(p_1-q_{\sigma(3)})-z(p_3-q_{\sigma(1)})|,
\end{array}
\right.
\end{equation}
where $s=k_1-k_2$ and $h$ has the sign in \eqref{e:h}.
Since $(a,b)$ and $(a+np+m\Id,b+nq+m\Id)$, $n,m\in\Z$, induce the same
circle action on $E_{p,q}$, we can assume $k_1=s$, $k_2=0$ and we get
\begin{equation}\label{e:ab}
\left\{\begin{array}{l}
a=\left(-z,-x-s(p_2-q_{\sigma(3)}),y+w+s(p_1-q_{\sigma(2)})\right),
\vspace{2ex} \\
b_{\sigma}=\left(w-x-s(p_2-q_{\sigma(3)}),y-z+s(p_1-q_{\sigma(2)}),0\right).
\end{array}
\right.
\end{equation}

Now, by \rref{r:afcond}, $\S^1_{a,b}$ is almost free if and only if the
orders in \eqref{e:orders} are nonzero. On the other hand, using
\eqref{free} it is easy to check that
$\kk_{\sigma\circ(12)},\kk_{\sigma\circ(23)},\kk_{\sigma\circ(13)}\neq0$,
unless at least two of the integers $p_i-q_j$ with $i\neq j$ are $\pm1$,
and $\kk_\sigma=0$ or $4$. Therefore, since $h$ is odd, we can make
$0,4\neq\kk_\sigma\leq h$ by choosing $s$ appropriately, which proves our
claim.
\end{proof}

\medskip

In order for an action to have only one singular point, there should be 3
vertices with the same parity that are regular, where recall that the
{\it parity} of the vertex $\hat\pi(C_\sigma)\in O$ is the parity of
$\sigma$. That is, there should exist 3 permutations $\sigma$, $\sigma'$,
$\sigma''$ with the same parity such that
$\kk_\sigma=\kk_{\sigma'}=\kk_{\sigma''}=1$. Notice that, by the last
observation in the proof of \pref{p:3lens}, in this situation the action
is automatically almost free. In addition, when such permutations exist,
the singular locus of $O$ is composed of at most the 3 vertices whose
parity is the opposite to that of $\sigma$.

To see when such an action exists, define the
following integers mod $h=|H^4(E_{p,q},\Z)|$:
\begin{equation}\label{e:alpha}
\alpha(\sigma,\e_1,\e_2):=\left((x+y-z)(p_1-q_{\sigma(1)})
-(w+z-x)(p_2-q_{\sigma(2)})+1\right)  \mod\, h
\end{equation}
where $x,y,z,w$ are defined in \eqref{e:cofactors}. One easily verifies
that $\alpha(\sigma,\e_1,\e_2)$ does not depend on the choice of
$x,y,z,w$, i.e., they do not depend on $k_1, k_2$, and are hence well
defined for every Eschenburg manifold $E_{p,q}$. Now, Theorem B in the
Introduction is an immediate consequence of the first equation in
\eqref{e:orders}.

Notice that the condition in Theorem B depends only on the parity of
$\sigma$, hence giving~8 tests to check. Moreover, it is trivially
satisfied when $h=1$, i.e., $H^4(E_{p,q},\Z)=0$. Using \eqref{e:h} one
easily sees that there are infinitely many Eschenburg manifolds with
$h=1$, for example, $p=(2k+2,k,0)$, $q=(2k+1,k+2,-1)$, with $k\in\Z$ such
that $\gcd(k-1,3)=1$. One can further try to minimize the singular set to
get only one point. It turns out that such actions exist in abundance.
For example, on $E_{(3,2,1),(4,2,0)}$, the almost free circle action
given by $a=(1,1,0)$, $b=(2,0,0)$ has only one singular point of order 3.

However, no Eschenburg manifold with positive curvature seems to satisfy
the condition in Theorem B, aside from the ones that already admit a free
circle action. By means of a computer program, we searched among all
positively curved spaces with $h\leq\ $100,000, in total 103,569,197
Eschenburg manifolds. It turns out that none of these spaces satisfies
the condition, apart from the 31,467 ones that admit free actions.
Moreover, none of these 31,467 spaces admit an isometric circle action
with only one singular point. On the other hand, we will see that there
are many spaces that admit circle actions with only two singular points
of opposite parity.

\smallskip

We can use the above methods to obtain a nice singular locus for a
general positively curved Eschenburg space. We give two typical
examples here, obtaining in particular the proof of Theorem C.

\smallskip

{\it Cohomogeneity one.}
Consider an arbitrary cohomogeneity one Eschenburg manifold, that is,
$E_d=E_{(1,1,d),(0,0,d+2)}$, $d\geq 0$. It has positive curvature if
$d>0$ and satisfies $h=2d+1$. For $d\leq2$, $E_d$ is known to admit a
free $\S^1$ action; cf. Section 2. So, assume $d\geq 3$. We want to study
all isometric $\ \S^1$ actions on $E_d$ for which at least 3 of the
vertices are regular. We will see that, for such an action, the regular
vertices do not have the same parity. As a consequence of this analysis,
we will obtain the proof of Theorem B in the introduction.

If there are 3 regular vertices, at least two of them correspond to
permutations $\sigma'$, $\sigma''$ with the same parity as the one
of, say, $\sigma$. Thus we are in the situation of the proof of
\pref{p:3lens}, and we have \eqref{e:cofactors}, \eqref{e:orders} and
\eqref{e:ab}. It is not difficult to compute all the possibilities in
\eqref{e:orders}, getting two cases:

{\it Case} $(a)$: $\sigma(3)=3$.
The integers $x=\e_1$, $z=-\e_2$, $y=w=0$ solve \eqref{e:cofactors},
we can assume that $\sigma=\Id$, and \eqref{e:orders} becomes
$$
\kk_{\Id}=|2\lambda-s|,\
\kk_{(12)}=|s|,\
\kk_{(23)}=|d\lambda-\e_2|, \
\kk_{(13)}=|d\lambda-\e_1|,
$$
with $\lambda=s(d+1)-\e_1-\e_2$. The only other order that can be 1 is
the second, $\kk_{(12)}$, that is, we can assume $s=1$ and then
$a =(\e_2,d+1-\e_1,1)$, $b=(d+1-\e_1,\e_2+1,0)$. Since there is no lens
space connecting $C_{(23)}$ and $C_{(13)}$, the only possibly singular
lens spaces are ${\mathcal L}_{11}$ and ${\mathcal L}_{22}$ whose orders
are $\gcd(d-1,3)$ and $\gcd(d-2,5)$ if $\e_1=-\e_2$, or both equal to 1
if $\e_1=\e_2$. Therefore, the minimal singular locus we can get in case
$(a)$ is given by~3 isolated singular points for $\e_1=\e_2=1$, whose
orbifold orders are $2d-3$, $d(d-1)-1$ and $d(d-1)-1$. Notice that for
$d\leq2$ we recover the known free $\ \S^1$ actions.

{\it Case} $(b)$: $\sigma(3)\neq 3$.
The integers $y=-\e_1$, $w=\e_2$, $x=z=0$ solve \eqref{e:cofactors},
we can assume that $\sigma(3)=1$, and we get for \eqref{e:orders}
$$
\kk_{\sigma}=|(2d+1)s-\e_1(d+1)+\e_2|,\
\kk_{\sigma\circ(12)}=|(d+1)(s-\e_1)+\e_2|,
$$
$$
\kk_{\sigma\circ(23)}=|2s-\e_1|,\
\kk_{\sigma\circ(13)}=|ds+\e_2|.
$$
Notice that, again, $\kk_\sigma > 1$, which says that no cohomogeneity
one Eschenburg manifold with $d\geq 3$ satisfies the condition in
Theorem B. To get one of the other 3 orders equal to one, we should have
either

\begin{itemize}
\item[$(b_1)$] $s=0$, with $a=(0,0,\e_2-\e_1)$,
$b_{\sigma}=(\e_2,-\e_1,0)$, and
$\kk_\sigma=\kk_{\sigma\circ(12)}=d+1-\e_1\e_2$ and the other 4
orders equal to one; or
\item[$(b_2)$] $s=\e_1$, with $a=(0,-\e_1,\e_2)$,
$b_{\sigma}=(\e_2-\e_1,0,0)$,
$\kk_\sigma=\kk_{\sigma\circ(23)}=d+\e_1\e_2$ and the other 4 orders
also equal to one.
\end{itemize}
For $(b_1)$, since $b_1=0$ and $b_2,b_3=\pm1$, the orbifold order of the
only possibly singular lens space ${\mathcal L}_{31}$ is
$\kk_{31}=\kk((1,-d-1,1),(-b_2,-b_3,-b_2))=d+1-\e_1\e_2=\kk_\sigma$.
So, it coincides with the orbifold order of the two exceptional circles
it contains. Therefore, the singular locus is a smooth totally geodesic
2-sphere with constant orbifold group $\Z_{d+1-\e_1\e_2}$.

Similarly, for $(b_2)$, the orbifold order of the only possible
singular lens space ${\mathcal L}_{1\sigma(1)}$ is
$\kk_{1\sigma(1)}=\kk((1,1-q_{\sigma(2)},d),(-\e_1,-\e_1,\e_2))
=\gcd(q_{\sigma(2)},d+\e_1\e_2)$.

For $\sigma=(13)$ we get as before that $\kk_{13}=d+\e_1\e_2=\kk_\sigma$.
So, also in this case the singular locus is a smooth totally geodesic
2-sphere with constant orbifold group $\Z_{d+\e_1\e_2}$, and the
exceptional set is the homogeneous lens space
${\mathcal L}_{13}=\Sph^3/\Z_{d+1}$. Notice also that this action gives
the smallest possible value for $\kk_\sigma$ under the assumption that
\mbox{$\kk_{\sigma'}=\kk_{\sigma''}=1$}.

On the other hand, for $\sigma=(123)$, it holds that
$\kk_{12}=\gcd(d+2,2-\e_1\e_2)$. Therefore, either $\kk_{12}=1$ if
$\e_1=\e_2$, in which case the singular locus has only two vertices
with the same orbifold group $\Z_{d+1}$, or $\kk_{12}=\gcd(d-1,3)$
if $\e_1=-\e_2$, for which we get a 2-sphere as singular locus if
and only if $3$ divides $d-1$, that is smooth if only if $d=4$.

\bigskip

{\it Cohomogeneity two.}
Observe that $(p,q)=((1,c,d),(0,0,c+d+1))$ lies in the same plane
generated in Theorem B $iv)$ for $E_d$. So, for any element in this
cohomogeneity two subfamily, the same action has as its exceptional set
the smooth lens space ${\mathcal L}_{13}=\Sph^3/\Z_{d+c}$ with constant
isotropy group $\Z_{d-c}$.

For a general cohomogeneity two Eschenburg manifold $E_{p,q}$, that is,
$p=(c,d,e)$, \mbox{$q=(0,0,c+d+e)$}, with
$\gcd(c,d)=\gcd(d,e)=\gcd(e,c)=1$, we can consider the action generated
by $a=(0,0,0)$, $b=(1,-1,0)$. In this case,
we get $\kk_{\Id}=\kk_{(12)}=|c+d|$, $\kk_{(123)}=\kk_{(23)}=|c+e|$,
$\kk_{(13)}=\kk_{(132)}=|d+e|$. The orders of the lens spaces are
$\kk_{13}=\gcd(2,d+e)$, $\kk_{23}=\gcd(2,c+e)$, $\kk_{33}=\gcd(2,c+d)$,
and $\kk_{11}=\kk_{12}=\gcd(c+d,c+e)$, $\kk_{21}=\kk_{22}=\gcd(c+d,d+e)$,
$\kk_{31}=\kk_{32}=\gcd(c+e,d+e)$. In particular, for the positively
curved Eschenburg space given by $(c,d,e)=(1,2,3)$ the singular locus
consists of one 2-sphere and 4 isolated points, whereas for
$(c,d,e)=(1,2,-3)$ it consists of one 2-sphere and 2 isolated points. On
the other hand, $(c,d,e)=(1,3,5)$ has the full set in Figure 1 as its
singular locus.

\begin{rem*}
With a simplified version of the arguments in Section 3, it is easy to
compute the singular locus of the circle action that defines a general
Eschenburg orbifold $E_{p,q}$. Again, it is given by Figure 1, composed
of the $C_\sigma$'s and the ${\mathcal L}_{ij}$'s, but now the orders of
their cyclic orbifold groups are $\gcd(p-q_\sigma)$ and
$\gcd(p-q_\sigma,p-q_{\sigma'})$, respectively. In particular, unlike in
the case of circle actions on Eschenburg spaces, the order of the lens
spaces is the gcd of the orders of the two circles which are contained in
them. There is now no difficulty to obtain positively curved Eschenburg
orbifolds that are not manifolds with the smallest possible singular
locus. For example, $E_{(5,3,-5),(2,1,0)}$ has only one circle as
singular locus with order~3.
\end{rem*}

\providecommand{\bysame}{\leavevmode\hbox to3em{\hrulefill}\thinspace}

\end{document}